\documentclass[a4paper, 12pt]{amsart}
\usepackage{amsmath, amssymb, eucal, amscd, amstext, enumerate}

\newtheorem{thm}{Theorem}%[section]
\newtheorem{lem}[thm]{Lemma}

\newtheorem{prop}[thm]{Proposition}
\newtheorem{cor}[thm]{Corollary}

\newenvironment{prf}{{\noindent \textbf{Proof:}\ }}{\hfill $\Box$\\ \smallskip}

\numberwithin{equation}{section}

\newcommand{\smnoind}{\smallskip\noindent}

\newcommand{\CE}{\mathcal{E}}
\newcommand{\CB}{\mathcal{B}}
\newcommand{\CH}{\mathcal{H}}

\begin{document}

\title[]{Property $T$ of reduced $C^*$-crossed products by discrete groups}

\author{Baojie Jiang \and Chi-Keung Ng}

\address[Baojie Jiang]{Shanghai Center for Mathematical Sciences, Fudan University, Shanghai 200433, China.}
\email{jiangbaojie@gmail.com}
\address[Chi-Keung Ng]{Chern Institute of Mathematics and LPMC, Nankai University, Tianjin 300071, China.}
\email{ckng@nankai.edu.cn; ckngmath@hotmail.com}

\thanks{The second named author is supported by the National Natural Science Foundation of China (11471168)}



\begin{abstract}
We generalize the main result of \cite{Kamalov} and show that if $G$ is an amenable discrete group with an action $\alpha$ on a finite nuclear unital $C^*$-algebra $A$ such that the reduced crossed product $A\rtimes_{\alpha,r} G$ has property $T$, then $G$ is finite and $A$ is finite dimensional.
As an application, an infinite discrete group $H$ is non-amenable if and only if the uniform Roe algebra $C^*_u(H)$ has property $T$.
\end{abstract}

\maketitle

\section{Introduction}

Property $T$ for unital $C^*$-algebras was introduced by Bekka in \cite{Bekka} and was studied by different people (see e.g. \cite{Brown,Kamalov,LN,Ng}).
In particular, it was shown by Kamalov in \cite{Kamalov} that
\begin{quotation}
if $G$ is a discrete amenable group acting on a commutative unital $C^*$-algebra $A$ such that the crossed product has property $T$, then $G$ is finite and $A$ is finite dimensional.
\end{quotation}
The aims of this paper is to extend this result to the case of finite nuclear unital $C^*$-algebras, and to give an application of this result.
As expected, a result of Brown in \cite{Brown} is one of our main tools.

\section{The main results}

Throughout this article,  $G$ is a discrete group acting on a unital $C^*$-algebra $A$ through an action $\alpha$ (by automorphisms).

Let $T(A)$ be the set of all traical states on $A$.
For any $\tau\in T(A)$, we denote by $\pi_\tau:A \to \CB(\CH_\tau)$ the GNS representation corresponding to $\tau$ and by $\xi_\tau$ a norm one cyclic vector in $\CH_\tau$ with
$$\tau(a) = \langle \pi_\tau(a)\xi_\tau, \xi_\tau \rangle \qquad (a\in A).$$
Recall that $A$ is said to be \emph{finite} if $T(A)$ separates points of $A_+$ (\cite[Theorem 3.4]{CP79}).
We also recall from \cite[Remark 2]{Bekka} that if $T(A) = \emptyset$, then $A$ has property $T$.

We use $T_\alpha(A)$ to denote the set of all $\alpha$-invariant traical states on $A$, and recall that $A$ is said to be
\emph{$\alpha$-finite} if $T_\alpha(A)$ separates points of $A_+$ (see \cite[Theorem 8.1]{CP79}).
We also denote by $A\rtimes_{\alpha,r} G$ the reduced crossed product of $\alpha$, and identify $A\subseteq A\rtimes_{\alpha,r} G$ as well as $G\subseteq A\rtimes_{\alpha,r} G$ through their canonical embeddings

Let us first give the following well-known facts.
Since we cannot find precise references for them, we present their simple arguments here.

\begin{lem}\label{lem:inv-st}
(a) $T(A\rtimes_{\alpha,r} G)\neq \emptyset$ if and only if $T_\alpha(A)\neq \emptyset$.

\smnoind
(b) If $A$ is $\alpha$-finite, then $A\rtimes_{\alpha,r} G$ is finite.

\smnoind
(c)	If $G$ is amenable and $T(A)\neq \emptyset$, then $T_\alpha(A)\neq \emptyset$.
\end{lem}
\begin{prf}
Let us denote $B:=A\rtimes_{\alpha,r} G$, and consider $\CE:B \to A$ to be the canonical conditional expectation (see e.g. \cite[Proposition 4.1.9]{BO}).

\smnoind
(a) If $\sigma\in T(B)$, then $\sigma (\alpha_t(a)) = \sigma(t a t^{-1}) = \sigma(a)$ ($a\in A;t\in G$), which means that $\sigma|_A\in T_\alpha(A)$.
Conversely, for any $\tau\in T_\alpha(A)$ and any $x=\sum_{s\in G} a_s s$ with $a_s =0 $ except for a finite number of $s$, one has
$$\tau(\CE(x^*x)) = \tau\Big(\sum_{r\in G} \alpha_{r^{-1}}(a_r^* a_r)\Big) = \tau\Big(\sum_{r\in G} a_ra_r^*\Big) = \tau(\CE(xx^*)).$$
Hence, $\tau\circ \CE$ belongs to $T(B)$, because it is continuous.

\smnoind
(b) As $\CE$ is faithful, we know that $B$ is a Hilbert $A$-module under the $A$-valued inner product
$$\langle x, y\rangle_A:= \CE(x^*y) \qquad (x,y\in B).$$
Moreover, for any $\tau \in T_\alpha(A)$, if $\pi_\tau^B$ is the canoincal representation of $B$ on the Hilbert space $B\otimes_{\pi_\tau} \CH_\tau$ (see e.g. \cite[Proposition 4.5]{Lance} for its definition; note that we identify a Hilbert $\mathbb{C}$-module with a Hilbert space by considering the conjugation of the inner product), then $(B\otimes_{\pi_\tau} \CH_\tau, \pi_\tau^B)$ coincides with $(\CH_{\tau\circ \CE}, \pi_{\tau\circ \CE})$ (observe that $1\otimes \xi_\tau$ is a cyclic vector for $\pi_\tau^B$ with the state defined by $1\otimes \xi_\tau$ being $\tau\circ \CE$).

Let $(\CH_0, \pi_0):= \bigoplus_{\tau\in T_\alpha(A)} (\CH_\tau, \pi_\tau)$.
Since $A$ is $\alpha$-finite, one knows that $\pi_0$ is faithful.
It is easy to verify that the representation $\pi_0^B$ of $B$ on $B\otimes_{\pi_0} \CH_0$ induced by $\pi_0$ is also faithful, and that $\pi_0^B$ coincides with $\bigoplus_{\tau\in T_\alpha(A)} \pi_\tau^B$.
Consequently, $\bigoplus_{\tau\in T_\alpha(A)} (\CH_{\tau\circ \CE}, \pi_{\tau\circ \CE})$ is faithful, which means that $\{\tau\circ \CE: \tau \in T_\alpha(A)\}$ (which is a subset of $T(B)$ by the argument of part (a)) separates points of $B_+$.

\smnoind
(c)	Note that $T(A)$ is a non-empty weak$^*$-compact convex subset of $A^*$ and $\alpha$ induces an action of $G$ on $T(A)$ by continuous affine maps.
Day's fixed point theorem (see \cite[Theorem 1]{Day}) produces a fixed point $\tau_0\in T(A)$ for this action.
Obviously, $\tau_0\in T_\alpha(A)$.
\end{prf}

We warn the readers that part (c) of the above is not true for non-unital $C^*$-algebras.

Our main theorem concerns with the situation when $A\rtimes_{\alpha,r} G$ is nuclear and has property $T$.
In this situation, \cite[Theorem 5.1]{Brown} tells us that $A\rtimes_{\alpha,r} G$ is a direct sum of a finite dimensional $C^*$-algebra and a nuclear $C^*$-algebra with no tracial state %after submission
(note that although all $C^*$-algebras in \cite{Brown} are assumed to be separable, \cite[Theorem 5.1]{Brown} is true in the non-separable case because one can use \cite[Theorem 6.2.7]{BO} to replace \cite[Theorem 4.2]{Brown}).
The following theorem implies that if $G$ is infinite, then we arrive at one of the extreme that the whole reduced crossed product has no tracial state.
This proposition, together with its proof, is a main ingredient in the argument for our main theorem.

\begin{prop}\label{prop:nuc+prop-T}
Let $G$ be an infinite discrete group acting on a unital $C^*$-algebra $A$ through an action $\alpha$.
If $A\rtimes_{\alpha,r} G$ is nuclear and has property $T$, then $T(A\rtimes_{\alpha,r} G) = \emptyset$.
\end{prop}
\begin{prf}
Let $I_\alpha := \bigcap_{\tau\in T_\alpha(A)} \ker\pi_\tau$ and $A_\alpha:=A/I_\alpha$.
Suppose on contrary that $T(A\rtimes_{\alpha,r} G) \neq \emptyset$.
Then $I_\alpha\neq A$ because of Lemma \ref{lem:inv-st}(a).
As $\ker \pi_\tau= \{x\in A: \tau(x^*x) = 0\}$ ($\tau \in T(A)$), we know that $I_\alpha$ is $\alpha$-invariant, and hence $\alpha$ produces an action $\beta$ of $G$ on $A_\alpha$.
Moreover, every element in $T_\alpha(A)$ induces an element in $T_\beta(A_\alpha)$, which gives the $\beta$-finiteness of $A_\alpha$.

Since $A_\alpha\rtimes_{\beta,r}G$ is a quotient $C^*$-algebra of $A\rtimes_{\alpha,r}G$, the hypothesis implies $A_\alpha\rtimes_{\beta,r} G$ to be nuclear and having property $T$.
Therefore, \cite[Theorem 5.1]{Brown} tells us that $A_\alpha\rtimes_{\beta,r} G = C\oplus D$, where $C$ is finite dimensional and $T(D) = \emptyset$.
However, the finiteness of $A_\alpha\rtimes_{\beta,r} G$ (which follows from Lemma \ref{lem:inv-st}(b)) tells us that $D=(0)$.
Consequently, $A_\alpha\rtimes_{\beta,r} G$ is a non-zero finite dimensional $C^*$-algebra, which contradicts the fact that $G$ is infinite.
\end{prf}

The following is our main theorem which concerns with the other extreme.
More precisely, what we obtained is a situation (which include the one in \cite{Kamalov}) under which the reduced crossed product is finite dimensional.

Notice that the finiteness assumption of $A$ is indispensable.
In fact, if $A$ is the direct sum of $\mathbb{C}$ with a nuclear unital $C^*$-algebra having no tracial state, then $A$ has a tracial state (but is not finite), and the reduced crossed product of the trivial action of a finite group on $A$ is nuclear and has property $T$.
We will see at the end of this article that one cannot weaken the amenability assumption of $G$ neither.

\begin{thm}\label{thm:amen-gp-act-fin-nucl}
Let $G$ be an amenable discrete group and $A$ be a finite nuclear unital $C^*$-algebra.
If there is an action $\alpha$ of $G$ on $A$ such that $A\rtimes_{\alpha,r} G$ has property $T$, then $G$ is finite and $A$ is finite dimensional.
\end{thm}
\begin{prf}
Set $I_\alpha := \bigcap_{\tau\in T_\alpha(A)} \ker\pi_\tau$ and $A_\alpha:=A/I_\alpha$. Denote $B:= A\rtimes_{\alpha,r} G$.
The finiteness assumption of $A$ and Lemma \ref{lem:inv-st}(c) imply that $I_\alpha\neq A$ and that $T(B)\neq \emptyset$ (see also Lemma \ref{lem:inv-st}(a)).
Hence, $G$ is finite (by Proposition \ref{prop:nuc+prop-T}).
Moreover, the argument of Proposition \ref{prop:nuc+prop-T} tells us that $I_\alpha$ is $\alpha$-invariant and $B_\alpha:= A_\alpha\rtimes_{\beta,r} G$ is finite dimensional.
Therefore, it suffices to show that $I_\alpha = \{0\}$.

Suppose on the contrary that $I_\alpha \neq \{0\}$.
By \cite[Theorem 5.1]{Brown}, we know that $B \cong B_0\oplus B_1$, where $B_0$ is finite dimensional and $T(B_1)=\emptyset$.
Thus, $I_\alpha\rtimes_{\alpha,r} G = J_0 \oplus J_1$, with $J_k$ being a closed ideal of $B_k$ for $k\in\{0,1\}$.
The short exact sequence
$$0\to I_\alpha \to A \to A_\alpha \to 0,$$
induces a short exact sequence concerning their full crossed products, which coincide with the reduced crossed products because $G$ is amenable.
From this, we obtain
$$B_\alpha = B/(I_\alpha\rtimes_{\alpha,r} G) = B_0/J_0 \oplus B_1/J_1.$$
Hence, $B_1/J_1$ is a quotient $C^*$-algebra of the finite dimensional $C^*$-algebra $B_\alpha$, which implies $J_1 = B_1$ (otherwise, $B_1$ will have a tracial state).
Consequently, $B_\alpha \cong B_0/J_0$, or equivalently, $B_0 \cong B_\alpha \oplus J_0$ (as $B_0$ is finite dimensional).
This gives
$$B \cong B_\alpha \oplus J_0 \oplus B_1 = B_\alpha \oplus (I_\alpha\rtimes_{\alpha,r} G).$$
Thus, $I_\alpha\rtimes_{\alpha,r} G$ is unital and so is $I_\alpha$ (but its identity may not be the identity of $A$).

Now, by the finiteness assumption of $A$, one knows that $T(I_\alpha)\neq \emptyset$, and Lemma \ref{lem:inv-st}(c) produces an element $\tau\in T_\alpha(I_\alpha)$.
Let $\Phi:A\to I_\alpha$ be the canonical $G$-equivariant $^*$-epimorphism, and define
$$\tau'(a):=  \langle \pi_\tau(\Phi(a)) \xi_\tau, \xi_\tau\rangle \qquad (a\in A).$$
Then $\tau'\in T_\alpha(A)$ and $\tau'|_{I_\alpha} = \tau$.
However, the existence of $\tau'$ contradicts the definition of $I_\alpha$.
\end{prf}

\begin{cor}\label{cor:non-amen<->prop-T}
Let $G$ be an infinite discrete group and $\alpha_G$ be the left translation action of $G$ on $\ell^\infty(G)$.
The following are equivalent.
\begin{enumerate}
\item $G$ is non-amenable.

\item $\ell^\infty(G)\rtimes_{\alpha_G,r} G$ does not have a tracial state.

\item $\ell^\infty(G)\rtimes_{\alpha_G,r} G$ has strong property $T$ (see \cite{LN}).

\item $\ell^\infty(G)\rtimes_{\alpha_G,r} G$ has property $T$.

%after submission
\item There is a finite nuclear unital $C^*$-algebra $A$ and an action $\alpha$ of $G$ on $A$ such that $A\rtimes_{\alpha, r}G$ has property $T$.
\end{enumerate}
\end{cor}
\begin{prf}
If $G$ is non-amenable, then $T_{\alpha_G}(\ell^\infty(G)) = \emptyset$ and Lemma \ref{lem:inv-st}(a) tells us that Statement (2) holds.
On the other hand, if $\ell^\infty(G)\rtimes_{\alpha_G,r} G$ does not have a tracial state, then \cite[Proposition 5.2]{LN} gives Statement (3).
Moreover, a strong property $T$ $C^*$-algebra clearly have property $T$.
Finally, suppose that %after submission
$A\rtimes_{\alpha,r} G$ has property $T$ but $G$ is amenable.
Then Theorem \ref{thm:amen-gp-act-fin-nucl} produces the contradiction that $G$ is finite.
\end{prf}

The following comparison of Corollary \ref{cor:non-amen<->prop-T} with the main result of Ozawa in \cite{Ozawa} (see also Theorem 5.1.6 and Proposition 5.1.3 of \cite{BO}) may be worth mentioning:
\begin{quotation}
a discrete $G$ is exact if and only if $\ell^\infty(G)\rtimes_{\alpha_G,r} G$ is nuclear (or equivalently, the action $\alpha_G$ is amenable).
\end{quotation}

This result tells us that one cannot weaken the amenability assumption of $G$ in Theorem \ref{thm:amen-gp-act-fin-nucl} to an amenable action $\alpha$ with $A\rtimes_{\alpha,r} G$ being nuclear, since if $G$ is an infinite exact non-amenable group, the action of $G$ on $\ell^\infty(G)$ is amenable, and the reduced crossed product has property $T$ and is nuclear.


\begin{thebibliography}{99}

\bibitem{Bekka}
B. Bekka, Property (T) for $C^* $-algebras, Bull. London Math. Soc. \textbf{38} (2006), 857-867.

\bibitem{Brown}
N.P. Brown, Kazhdan's property T and $C^\ast$-algebras, J. Funct. Anal., \textbf{240} (2006), 290-296.

\bibitem{BO}
N.P. Brown and N. Ozawa, \emph{$C^* $-algebras and finite dimensional approximations}, Grad. Stud. Math. \textbf{88}, American Mathematical Society, Providence, RI, (2008).

\bibitem{CP79}
J.\ Cuntz and G.\ K.\ Pedersen,
Equivalence and traces on $C^*$-algebras, J.\ Funct.\ Anal., \textbf{33} (1979), 135-164.

\bibitem{Day}
M.M. Day, Fixed-point theorems for compact convex sets, Illinois J. Math. \textbf{5} (1961), 585-590.

\bibitem{Kamalov}
F. Kamalov, Property $T$ and Amenable Transformation Group $C^*$-algebras, Canad. Math. Bull. \textbf{58} (2015), 110-114.

\bibitem{Lance}
E.C. Lance, \emph{Hilbert $C^*$-modules, A toolkit for operator algebraists}, London Math. Soc. Lect. Note Ser. \textbf{210}, Camb. Univ. Press (1995).

%\bibitem{LNW}
%C.W. Leung, C.K. Ng and N.C. Wong, Property $(T)$ for non-unital $C^*$-algebras, J. Math. Anal. Appl. \textbf{341} (2008), 1102–1106.

\bibitem{LN}
C.W. Leung and C.K. Ng, Property $(T)$ and strong property $(T)$ for unital $C^*$-algebras, J. Funct. Anal. \textbf{256} (2009), 3055–3070.

\bibitem{Ng}
C.K. Ng, Property $T$ for general $C^*$-algebras, Math. Proc. Camb. Philos. Soc. \textbf{156} (2014), 229–239.

\bibitem{Ozawa}
N. Ozawa, Amenable actions and exactness for discrete groups,  C. R. Acad. Sci. Paris S\'{e}r. I Math. \textbf{330} (2000), 691–695.

\end{thebibliography}
\end{document}